\documentclass[a4paper,10pt]{article}
\usepackage{stmaryrd}
\usepackage{amsfonts}
\usepackage{bbm}
\usepackage{amscd}
\usepackage{mathrsfs}
\usepackage{latexsym,amssymb,amsmath,amscd,amscd,amsthm,amsxtra,xypic}
\usepackage[dvips]{graphicx}
\usepackage[utf8]{inputenc}
\usepackage[T1]{fontenc}
\usepackage{lmodern}
\usepackage{amssymb}
\usepackage[all]{xy}
\usepackage{nicefrac,mathtools,enumitem}
\usepackage{microtype}

\newtheorem{thm}{Theorem}[section]
\newtheorem{lem}[thm]{Lemma}
\newtheorem{cor}[thm]{Corollary}
\newtheorem{pro}[thm]{Proposition}
\newtheorem{ex}[thm]{Example}
\newtheorem{rmk}[thm]{Remark}
\newtheorem{defi}[thm]{Definition}

\newcommand{\be }{\begin{equation}}
\newcommand{\ee }{\end{equation}}

\newcommand{\pf}{\noindent{\bf Proof.}\ }
\newcommand {\emptycomment}[1]{} 
\newcommand{\g}{\frkg}

\newcommand{\huaG}{\mathcal{G}}
\newcommand{\adj}{\frka\frkd}

\newcommand{\frka}{\mathfrak a}

\newcommand{\frkd}{\mathfrak d}

\newcommand{\frkf}{\mathfrak f}
\newcommand{\frkg}{\mathfrak g}

\newcommand{\frki}{\mathfrak i}

\def\qed{\hfill ~\vrule height6pt width6pt depth0pt}

\newcommand{\half}{\frac{1}{2}}

\newcommand{\Dorfman}[1]{\left\llbracket  #1\right\rrbracket }
\newcommand{\Courant}[1]{\{ #1\}}

\newcommand{\id}{\mathrm{id}}

\newcommand{\br}[1]{   [ \cdot,    \cdot  ]   }

\newcommand{\Hom}{\mathrm{Hom}}

\newcommand{\gl}{\mathfrak {gl}}
\newcommand{\ol}{\mathfrak {ol}}

\newcommand{\ad}{\mathrm{ad}}
\newcommand{\pr}{\mathrm{pr}}

\newcommand{\img}{\mathrm{im}}

\newcommand{\sgn}{\mathrm{sgn}}

\begin{document}
\title{
{Remarks on   Leibniz algebras
\thanks
 {
Research partially supported by  NSF of China (11101179).
 }
} }
\author{Yunhe Sheng  \\
Department of Mathematics, Jilin University,
 Changchun 130012,  China
\\\vspace{3mm}
email: shengyh@jlu.edu.cn\\
Zhangju Liu\\
Department of Mathematics, Peking University,  Beijing 100871,
China\\\vspace{3mm} email: liuzj@pku.edu.cn}
\date{}
\maketitle

\begin{abstract}
In this paper, first  we construct a Lie 2-algebra  associated to
every Leibniz algebra via the skew-symmetrization. Furthermore, we
introduce the notion of the naive representation for a Leibniz
algebra in order to realize the abstract operations as a concrete
linear operation. At last, we study some properties of naive
cohomologies.
\end{abstract}


\section{The skew-symmetrization of a Leibniz algebra}

In this section, we construct  a Lie 2-algebra from a Leibniz
algebra via the skew-symmetrization. The notion of Leibniz algebras
was introduced by Loday in \cite{Loday}, which  is a vector space
$\frkg$, endowed with a linear map
$[\cdot,\cdot]_\frkg:\frkg\otimes\frkg\longrightarrow\frkg$
satisfying
\begin{equation}\label{eq:Leibniz}
[x,[y,z]_\frkg]_\frkg=[[x,y]_\frkg,z]_\frkg+[y,[x,z]_\frkg]_\frkg,\quad
\forall~x,y,z\in \frkg.
\end{equation}
The {\bf left center} is given by
\begin{equation}
  \mathrm Z(\g)=\{x\in\g|~ [x,y]_\g=0,\quad\forall y\in\g\}.
\end{equation}
It is obvious that $\mathrm Z(\g)$ is an ideal and the quotient
Leibniz algebra $\g/ \mathrm Z(\g)$ is actually a Lie algebra since
$[x,x]_\g\in\mathrm Z(\g)$, for all $x \in\g$.

A Lie 2-algebra is a categorification of a Lie algebra, which is
equivalent to a 2-term $L_\infty$-algebra (see
\cite{baez:2algebras}, \cite{Stasheff} for more details).

\begin{defi}\label{defi:Lie 2}
A {\bf Lie $2$-algebra} is a graded vector space $\g=\g_{1}
\oplus\g_{0}$, together  with linear maps
$\{l_k:\wedge^k\frkg\longrightarrow \frkg,k=1,2,3\}$ of degrees $\deg(l_k)=k-2$
satisfying the following equalities:
\begin{enumerate}
\item[\rm(a)] $l_1l_2(x,a)=l_2(x,l_1( a))$,
\item[\rm(b)] $l_2(l_1( a),b)=l_2(a,l_1( b))$,
\item[\rm(c)] $l_2(x,l_2(y,z))+c.p.=l_1 l_3(x,y,z)$,
\item[\rm(d)] $l_2(x,l_2(y,a))+l_2(y,l_2(a,x))+l_2(a,l_2(x,y))=l_3(x,y,l_1( a))$,
\item[\rm(e)] $l_3(l_2(x,y),z,w)+c.p.=l_2(l_3(x,y,z),w)+c.p.$,
\end{enumerate}
for all $ x,y,z,w\in \g_0,~a,b\in \g_{1}$, where  $c.p.$ means cyclic
permutations.
\end{defi}

Given  a Leibniz algebra  $(\g,[\cdot,\cdot]_\g)$, introduce the
following skew-symmetric bracket on $\g$:
\begin{equation}
\Dorfman{x,y}=\half([x,y]_\g-[y,x]_\g),\quad\forall x,y\in\g,
\end{equation}
and denote by $J_{x,y,z}$ the corresponding Jacobiator, i.e.
\begin{equation}\label{eq:skew}
  J_{x,y,z} =\Dorfman{x,\Dorfman{y,z}}+\Dorfman{y,\Dorfman{z,x}}+\Dorfman{z,\Dorfman{x,y}}.
\end{equation}

\begin{pro}\label{pro:J}
  Let $(\g,[\cdot,\cdot]_\g)$ be a Leibniz algebra.
  \begin{itemize}
    \item[\rm(i)] For all $x,y,z\in\g$, we have
  \begin{equation}
   J_{x,y,z}=\frac{1}{4}\big([[z,y]_\g,x]_\g+[[x,z]_\g,y]_\g+[[y,x]_\g,z]_\g\big).
  \end{equation}

     \item[\rm(ii)] $J_{x,y,z}\in \mathrm Z(\g)$, i.e. $J_{x,y,z}$ is in the left center of  $(\g,[\cdot,\cdot]_\g)$.

      \item[\rm(iii)] For all $x,y,z,w\in\g$, we have \begin{eqnarray}
       &&\nonumber \Dorfman{x,J_{y,z,w}}-\Dorfman{y,J_{x,z,w}}+\Dorfman{z,J_{x,y,w}}-\Dorfman{w,J_{x,y,z}}\\\label{eq:closed}
    &&-J_{\Dorfman{x,y},z,w}+J_{\Dorfman{x,z},y,w}-J_{\Dorfman{x,w},y,z}-J_{\Dorfman{ y,z},x,w}+J_{\Dorfman{ y,w},x,z}-J_{\Dorfman{z,w},x,y}=0.
      \end{eqnarray}
  \end{itemize}
 \end{pro}
 \pf The first conclusion is obtained by straightforward computations.
  For any $w\in\g$, by \eqref{eq:Leibniz} and the fact that for all $x \in\g
, [x,x]_\g\in\mathrm Z(\g) $, we have
\begin{eqnarray*}
 ~[ J_{x,y,z},w]_\g&=&\frac{1}{4}\big([[[z,y]_\g,x]_\g+[[x,z]_\g,y]_\g+[[y,x]_\g,z]_\g,w]_\g\big)\\
 &=&\frac{1}{4}\big([[z,[y,x]_\g]_\g-[y,[z,x]_\g]_\g-[[z,x]_\g,y]_\g+[[y,x]_\g,z]_\g,w]_\g\big)\\
 &=&0,
\end{eqnarray*}
which implies that $J_{x,y,z}\in \mathrm Z(\g)$. At last, since the
bracket $\Dorfman{\cdot,\cdot}$ given by \eqref{eq:skew} is
skew-symmetric,  we have
\begin{eqnarray*}
 &&\Dorfman{x,J_{y,z,w}}-\Dorfman{y,J_{x,z,w}}+\Dorfman{z,J_{x,y,w}}-\Dorfman{w,J_{x,y,z}}\\
    &&-J_{\Dorfman{x,y},z,w}+J_{\Dorfman{x,z},y,w}-J_{\Dorfman{x,w},y,z}-J_{\Dorfman{ y,z},x,w}+J_{\Dorfman{ y,w},x,z}-J_{\Dorfman{z,w},x,y}\\
    &=&\underline{\Dorfman{x,J_{y,z,w}}}\underbrace{-\Dorfman{y,J_{x,z,w}}}+\underline{\underline{\Dorfman{z,J_{x,y,w}}}}
    -\underline{\underline{\underline{\Dorfman{w,J_{x,y,z}}}}}
    \end{eqnarray*}

    \begin{eqnarray*}
    &&-\Dorfman{\Dorfman{x,y},\Dorfman{z,w}}-\underline{\underline{\Dorfman{z,\Dorfman{w,\Dorfman{x,y}}}}}
    -\underline{\underline{\underline{\Dorfman{w,\Dorfman{\Dorfman{x,y},z}}}}}\\
    &&+\Dorfman{\Dorfman{x,z},\Dorfman{y,w}}
    +\underbrace{\Dorfman{y,\Dorfman{w,\Dorfman{x,z}}}}+\underline{\underline{\underline{\Dorfman{w,\Dorfman{\Dorfman{x,z},y}}}}}\\
    &&-\Dorfman{\Dorfman{x,w},\Dorfman{y,z}}\underbrace{-\Dorfman{y,\Dorfman{z,\Dorfman{x,w}}}}-\underline{\underline{\Dorfman{z,\Dorfman{\Dorfman{x,w},y}}}}
    \\&&-\Dorfman{\Dorfman{y,z},\Dorfman{x,w}}\underline{-\Dorfman{x,\Dorfman{w,\Dorfman{y,z}}}}-\underline{\underline{\underline{\Dorfman{w,\Dorfman{\Dorfman{y,z},x}}}}}\\
    &&+\Dorfman{\Dorfman{y,w},\Dorfman{x,z}}+\underline{\Dorfman{x,\Dorfman{z,\Dorfman{y,w}}}}+\underline{\underline{\Dorfman{z,\Dorfman{\Dorfman{y,w},x}}}}
    \\&&-\Dorfman{\Dorfman{z,w},\Dorfman{x,y}}\underline{-\Dorfman{x,\Dorfman{y,\Dorfman{z,w}}}}\underbrace{-\Dorfman{y,\Dorfman{\Dorfman{z,w},x}}}\\
    &=&-\Dorfman{\Dorfman{x,y},\Dorfman{z,w}}+\Dorfman{\Dorfman{x,z},\Dorfman{y,w}}-\Dorfman{\Dorfman{x,w},\Dorfman{y,z}}
    \\&&-\Dorfman{\Dorfman{y,z},\Dorfman{x,w}}+\Dorfman{\Dorfman{y,w},\Dorfman{x,z}}-\Dorfman{\Dorfman{z,w},\Dorfman{x,y}}\\
    &=&0.
\end{eqnarray*}
The proof is finished.\qed

\emptycomment{\begin{cor}
 With the above notations, $J$ is totally skew-symmetric.
\end{cor}
\pf For all $x,z\in\g$, by \eqref{eq:center}, we have
\begin{eqnarray*}
  J_{x,x,z}=\frac{1}{4}\big([[z,x]_\g,x]_\g+[[x,z]_\g,x]_\g+[[x,x]_\g,z]_\g\big)=0.
\end{eqnarray*}
Similarly, we have $J_{x,z,x}=0$ and $J_{z,x,x}=0$. Thus, $J$ is totally skew-symmetric.
\qed}

Next, for  a Leibniz algebra $(\g,[\cdot,\cdot]_\g)$,  we consider
the graded vector space $\huaG=\mathrm Z(\g)\oplus \g$, where
$\mathrm Z(\g)$ is of degree 1, $\g$ is of degree 0. Define a degree
$-1$ differential $l_1=\frki: \mathrm Z(\g)\longrightarrow \g$, the
inclusion. Define a degree $0$ skew-symmetric bilinear map $l_2$ and
a degree $1$ totally skew-symmetric trilinear map $l_3$ on $\huaG$
by

\begin{equation}
 \left\{ \begin{array}{ll}l_2(x,y)=\Dorfman{x,y}=\half([x,y]_\g-[y,x]_\g)&\forall~ x,y\in\g, \\
 l_2(x,c)=-l_2(c,x)=\Dorfman{x,c}=\half[x,c]_\g&\forall~ x\in\g,~ c\in\mathrm Z(\g), \\
  l_2(c_1,c_2)=0&\forall~  c_1,c_2\in\mathrm Z(\g), \\
  l_3(x,y,z)=J_{x,y,z}&\forall~ x,y,z\in\g.
 \end{array}\right.
\end{equation}

The following theorem is  our main result in this section, which says that one can obtain a Lie 2-algebra via the skew-symmetrization of a Leibniz algebra.
\begin{thm}\label{thm:main1}
  With the above notations, $(\huaG,l_1,l_2,l_3)$ is a Lie $2$-algebra.
\end{thm}
\pf By the definition of $l_1,~l_2$ and $l_3$, it is obvious that Conditions (a)-(d) in Definition \ref{defi:Lie 2} hold. By (iii) in Proposition \ref{pro:J}, Condition (e) also holds. Thus, $(\huaG,l_1,l_2,l_3)$ is a Lie $2$-algebra.\qed

\section{Representations of Leibniz algebras}
The theory of  representations and cohomologies of Leibniz algebras
was introduced and studied  in \cite{Loday and Pirashvili}.
Especially, faithful representations and conformal representations
of Leibniz algebras were studied in \cite{faithful rep} and
\cite{conformal rep} respectively. See   \cite{FMM,HuPeiLiu} for
more applications of cohomologies of Leibniz algebras.

\begin{defi} \label{defi:rep old}
A representation of the Leibniz algebra
$(\frkg,[\cdot,\cdot]_\frkg)$ is a triple  $(V,l,r)$, where  $V$ is
a vector space equipped with two linear maps
$l:\g\longrightarrow\gl(V)$ and $r:\g\longrightarrow\gl(V)$ such
that the following equalities hold:
\begin{equation}\label{condition of rep}
l_{[x,y]_\frkg}=[l_{x},l_{y}],\quad
r_{[x,y]_\frkg}=[l_{x},r_{y}],\quad r_{y}\circ l_{x}=-r_{y}\circ
r_{x}, ~~~ \, ~~~\, ~~~  \forall x,y\in\frkg.
\end{equation}
\end{defi}

\begin{defi} \label{defi:cohomo old}
The Leibniz cohomology of $\frkg$ with coefficients in $V$ is the
cohomology of the cochain complex
$C^k(\frkg,V)=\Hom(\otimes^k\frkg,V), (k\geq0)$ with the coboundary
operator
 $$\partial:C^k(\frkg,V)\longrightarrow C^{k+1}(\frkg,V)$$
defined by
\begin{eqnarray}
\nonumber\partial
c^k(x_1,\dots,x_{k+1})&=&\sum_{i=1}^k(-1)^{i+1}l_{x_i}(c^k(x_1,\dots,\widehat{x_i},\dots,x_{k+1}))
+(-1)^{k+1}r_{x_{k+1}}(c^k(x_1,\dots,x_k))\\
\label{formulapartial}&&+\sum_{1\leq i<j\leq
k+1}(-1)^ic^k(x_1,\dots,\widehat{x_i},\dots,x_{j-1},[x_i,x_j]_\frkg,x_{j+1},\dots,x_{k+1}).
\end{eqnarray}
The resulting cohomology is denoted by $H^\bullet(\g;l,r)$.
\end{defi}

Obviously, $(\mathbb R, 0, 0)$ is a representation of $\g$, which is
called the trivial representation. Denote the resulting cohomology
by $H^\bullet(\g).$ Another important representation is the adjoint
representation  $(\g, \ad_L, \ad_R)$, where
$\ad_L:\g\longrightarrow\gl(\g)$ and
$\ad_R:\g\longrightarrow\gl(\g)$  are defined as follows:
\begin{equation}
  \ad_L(x)(y)=[x,y]_\g,\quad \ad_R(x)(y)=[y,x]_\g, ~~~ \, ~~~\, ~~~  \forall x,y\in\frkg.
\end{equation}
 The resulting cohomology is denoted by
 $H^\bullet(\g;\ad_L, \ad_R).$\vspace{3mm}

The graded vector space $\bigoplus_k C^k(\frkg,\g)$ equipped with
the graded bracket
\begin{equation}
  [\alpha,\beta]=\alpha\circ\beta+(-1)^{pq+1}\beta\circ\alpha, \quad\forall ~\alpha\in C^{p+1}(\frkg,\g),~\beta\in C^{q+1}(\frkg,\g)
\end{equation}
  is a   graded Lie algebra, where $\alpha\circ\beta\in C^{p+q+1}(\g,\g)$ is defined by
\begin{eqnarray*}
  \alpha\circ\beta(x_1,\cdots,x_{p+q+1})&=&\sum_{k=0}^{p}(-1)^{kq}\Big(\sum_{\sigma\in sh(k,q)}\sgn(\sigma)\alpha(x_{\sigma(1)},\cdots,x_{\sigma(k)},\\
  &&\qquad\beta(x_{\sigma(k+1)},\cdots,x_{\sigma(k+q)},x_{k+q+1}),x_{k+q+2},\cdots,x_{p+q+1})\Big).
\end{eqnarray*}
See \cite{Balavoine,Fialowski} for more details. In particular, for
$\alpha\in C^2(\g,\g)$, we have
\begin{equation}
 ~[\alpha,\alpha](x,y,z)=2\alpha\circ\alpha(x,y,z)=2\Big(\alpha(\alpha(x,y),z)-\alpha(x,\alpha(y,z))+\alpha(y,\alpha(x,z))\Big).
\end{equation}
Thus, $\alpha$ defines a Leibniz algebra structure if and only if $[\alpha,\alpha]=0.$

For a representation  $(V,l,r)$ of $\g$, it is obvious that
$(V,l,0)$ is  a representation of $\g$ too. Thus,  we have two
semidirect product Leibniz algebras 
 $\g\ltimes_{(l,r)}V$  and $\g\ltimes_{(l,0)}V$ with the brackets $[\cdot,\cdot]_{(l,r)}$ and
$[\cdot,\cdot]_{(l,0)}$ respectively:
\begin{eqnarray*}
 ~ [x+u,y+v]_{(l,r)}&=&[x,y]_\g+l_xv+r_yu,\\
  ~ [x+u,y+v]_{(l,0)}&=&[x,y]_\g+l_xv.
\end{eqnarray*}
The right action $r$ induces a linear map $\overline{r}:(\g\oplus
V)\otimes (\g\oplus V)\longrightarrow \g\oplus V$ as follows:
$$
 \overline{r}(x+u,y+v)=r_yu.
$$

The adjoint actions maps $\ad_L$ and $\ad_R$ on the Leibniz algebra
$\g\ltimes_{(l,0)}V$ are given by
  $$\ad_L(x+u)(y+v)=[x,y]_\g+l_xv, ~~~ \, ~~~
  \ad_R(x+u)(y+v)=[y,x]_\g+l_yu.$$
\begin{pro}
  With the above notations, $\overline{r}$ satisfies the following Maurer-Cartan   equation on the Leibniz algebra $\g\ltimes_{(l,0)}V$:
  $$
  \partial \overline{r}-\half[\overline{r},\overline{r}]=0,
  $$
   where $\partial$ is  the coboundary
operator for the adjoint representation of $\g\ltimes_{(l,0)}V$.
 Consequently, the Leibniz algebra $\g\ltimes_{(l,r)}V$ is a deformation of the Leibniz algebra $\g\ltimes_{(l,0)}V$
  via the  Maurer-Cartan element $\overline{r}$.
  \end{pro}
\pf By direct computation, we have
\begin{eqnarray*}
\partial \overline{r}(x+u,y+v,z+w)&=&\ad_L(x+u)\overline{r}(y+v,z+w)-\ad_L(y+v)\overline{r}(x+u,z+w)\\
&&-\ad_R(z+w)\overline{r}(x+u,y+v)- \overline{r}([x+u,y+v]_{(l,0)},z+w)\\
 &&+\overline{r}(x+u,[y+v,z+w]_{(l,0)})-\overline{r}(y+v,[ x+u,z+w]_{(l,0)})\\
  &=&l_xr_zv-l_yr_zu-r_zl_xv+r_{[y,z]_\g}u-r_{[x,z]_\g}v.
\end{eqnarray*}
On the other hand, we have
\begin{eqnarray*}
  [\overline{r},\overline{r}](x+u,y+v,z+w)&=&2\big(\overline{r}(\overline{r}(x+u,y+v),z+w)-\overline{r}(x+u,\overline{r}(y+v,z+w))\\&&+\overline{r}(y+v,\overline{r}(x+u,z+w))\big)
  \\&=&2r_zr_yu.
\end{eqnarray*}
Thus, we have
\begin{eqnarray*}
  \Big(\partial \overline{r}-\half[\overline{r},\overline{r}]\Big)(x+u,y+v,z+w)&=&l_xr_zv-l_yr_zu-r_zl_xv+r_{[y,z]_\g}u-r_{[x,z]_\g}v-r_zr_yu\\
  &=&l_xr_zv-l_yr_zu-r_zl_xv+r_{[y,z]_\g}u-r_{[x,z]_\g}v+r_zl_yu\\
  &=&0.
\end{eqnarray*}
The proof is finished. \qed \vspace{3mm}

 Let  $(V^*,l^*,0)$ be the dual representation of the representation $(V,l,0)$ for
 $\g$,
  where $l^*:\g\longrightarrow\gl(V^*)$ is given by
  $$\langle l^*_x(\xi),u\rangle=-\langle\xi,l_xu\rangle,  ~~~~~ \, ~~~ \forall x \in \g, ~ \xi \in V^*, ~ u\in V.$$
  It is straightforward to see that $(V^*\otimes V,l^*\otimes 1+ 1\otimes l,0)$ is also a representation of $\g$,
  where $l^*\otimes 1+ 1\otimes l:\g\longrightarrow\gl(V^*\otimes V)$ is given by
$$
(l^*\otimes 1+ 1\otimes l)_x(\xi\otimes u)=(l^*_x\xi)\otimes u+\xi\otimes l_xu,\quad \forall \xi\in V^*,~u\in V.
$$

Since  $V^*\otimes V \cong \gl(V)$, an element $\xi\otimes u$ in
$V^*\otimes V$ can be identified with a linear map $A\in\gl(V)$ via
$A(v)=\langle\xi,v\rangle u$.

\begin{pro}\label{pro:1cocycle}
  With the above notations,  for all $A\in V^*\otimes V \cong \gl(V)$, we have
  $$
  (l^*\otimes 1+ 1\otimes l)_x A=[l_x,A].
  $$
  Moreover, the right action  $r:\g\longrightarrow \gl(V)$ is a $1$-cocycle for  the
   representation $(V^*\otimes V,l^*\otimes 1+ 1\otimes l,0)$ of $\g$.
\end{pro}
\pf Write $A=\xi\otimes u$, then we have
\begin{eqnarray*}
  (l^*\otimes 1+ 1\otimes l)_x A(v)&=&(l^*\otimes 1+ 1\otimes l)_x (\xi\otimes u)(v)=\Big((l^*_x\xi)\otimes u+\xi\otimes l_xu\Big)(v)\\
  &=&\langle l^*_x\xi,v\rangle u+\langle\xi,v\rangle l_xu=-\langle\xi,l_xv\rangle u+\langle\xi,v\rangle l_xu\\
  &=&[l_x,A](v).
\end{eqnarray*}
Therefore, we have
\begin{eqnarray*}
  \partial r(x,y)= (l^*\otimes 1+ 1\otimes l)_xr(y)-r([x,y]_\g)=[l_x,r_y]-r_{[x,y]_\g}=0,
\end{eqnarray*}
which implies that $r$ is a $1$-cocycle.\qed

\section{Naive representations of Leibniz algebras}

It is known that the aim of a representation is to  realize an
abstract algebraic structure  as a class of linear transformations
on a vector space.  Such as  a Lie algebra  representation is a
homomorphism from $\g$ to the general linear Lie algebra $\gl(V)$.
Unfortunately, the  representation of a Leibniz algebra discussed
above does not realize the abstract operation as a concrete linear
operation. Therefore, it is reasonable for us to  provide an
alternative definition for the representation of Leibniz algebras.
It is lucky that there is a god-given Leibniz algebra  worked as a
``general linear algebra''  defined as follows:

Given a vector space $V$,   then $(V,l=\id,r=0)$ is a natural
representation of $\gl(V)$, which is viewed as a Leibniz algebra.
 The corresponding semidirect product
Leibniz algebra structure on $\gl(V)\oplus V$ is given by
\begin{eqnarray*}
\Courant{A+u,B+v}&=&[A,B]+Av,\quad \forall~A,B\in \gl(V),~u,v\in V.
\end{eqnarray*}

 This Leibniz algebra is called {\bf
omni-Lie algebra} and denoted by $\ol(V)$.
 The notion of an omni-Lie algebra was introduced  by Weinstein in \cite{Weinomni}
 as the linearization of a Courant algebroid. The notion of a Courant algebroid was introduced
 in \cite{LWXmani}, which has been widely  applied in many
fields both for mathematics and physics (see \cite{Gualtieri, KS} for
more details). Its Leibniz algebra structure also  played  an
important role when studying the
 integrability of Courant brackets \cite{kinyon-weinstein}.

 Notice that  the skew-symmetric
 bracket,
 \begin{equation}
   \Dorfman{A+u,B+v}=[A,B]+\frac{1}{2}(Av-Bu),
 \end{equation}
which is obtained via the skew-symmetrization of $\Courant{\cdot,\cdot}$, is used
 in his original
 definition. As a special case in Theorem \ref{thm:main1},  $(\gl(V)\oplus
 V,\Dorfman{\cdot,\cdot})$ is a Lie 2-algebra.
Even though an  omni-Lie algebra is  not a Lie algebra,
  all Lie algebra structures on $V$ can be characterized by the Dirac structures
  in   $\ol(V)$. In fact, the next proposition will show that  every Leibniz algebra structure on $V$ can
  be realized as a Leibniz subalgebra of $\ol(V)$.
For any $\varphi:V\longrightarrow\gl(V)$, consider its graph
$$\huaG_\varphi=\{\varphi(u)+u\in\gl(V)\oplus V|~\forall u\in V\}.$$

\begin{pro}
  With the above notations, $\huaG_\varphi$ is a Leibniz subalgebra  of $\ol(V)$ if and only if
  \begin{equation}\label{eq:gr}
    [\varphi(u),\varphi(v)]=\varphi(\varphi(u)v),\quad\forall~u,v\in V.
  \end{equation}
  Furthermore, under this condition, $(V,[\cdot,\cdot]_\varphi)$ is a Leibniz algebra, where the linear map $[\cdot,\cdot]_\varphi:V\otimes V\longrightarrow V$ is given by
  \begin{equation}
    [u,v]_\varphi=\varphi(u)v,\quad \forall~u,v\in V.
  \end{equation}
\end{pro}
\pf Since $\ol(V)$ is a Leibniz algebra, we only need to show that $\huaG_\varphi$ is closed if and only if \eqref{eq:gr} holds. The conclusion follows from
\begin{eqnarray*}
  \Courant{\varphi(u)+u,\varphi(v)+v}=[\varphi(u),\varphi(v)]+\varphi(u)v.
\end{eqnarray*}
The other conclusion is straightforward. The proof is finished. \qed\vspace{3mm}

Recall that a representation of a Lie algebra $\g$ on a vector space
$V$ is a Lie algebra homomorphism from $\g$ to the Lie algebra
$\gl(V)$, which realizes an abstract Lie algebra as a subalgebra of
a concrete Lie algebra. Similarly, For a  Leibniz algebra, we
suggest the following definition:
\begin{defi}
  A naive representation of a Leibniz algebra $\g$ on a vector space $V$ is a Leibniz algebra homomorphism $\rho:\g\longrightarrow\ol(V)$.
\end{defi}

According to the two components of   $\gl(V)\oplus V$,  every linear
map $\rho:\g\longrightarrow\ol(V)$ can be split into two linear
maps: $\phi:\g\longrightarrow \gl(V)$ and
   $\theta:\g\longrightarrow V$. Then, we have


\begin{pro}\label{pro:rep con}
 A linear map $\rho=\phi+\theta:\g\longrightarrow\gl(V)\oplus V$ is a  naive
 representation  of $\g$
 if and only if $(V,\phi,0) $ is a representation of $\g$ and $\theta:\g\longrightarrow V$ is a
 $1$-cocycle for
the representation $(V,\phi,0)$.
\end{pro}
\pf On one hand, we have
$$
\rho([x,y]_\g)=\phi([x,y]_\g)+\theta([x,y]_\g).
$$
On the other hand, we have
\begin{eqnarray*}
  \Courant{\rho(x),\rho(y)}= \Courant{\phi(x)+\theta(x),\phi(y)+\theta(y)}=[\phi(x),\phi(y)]+\phi(x)\theta(y).
\end{eqnarray*}
Thus, $\rho$ is a homomorphism if and only if
\begin{eqnarray}\label{eq:rep con1}
   \phi([x,y]_\g)&=&[\phi(x),\phi(y)],\\
   \label{eq:rep con2}\theta([x,y]_\g)&=&\phi(x)\theta(y).
 \end{eqnarray}
By definition, Equalities \eqref{eq:rep con1}  and \eqref{eq:rep
con2} are equivalent to that $(V,\phi,0) $ is a representation and
$\theta:\g\longrightarrow V$ is a $1$-cocycle respectively. \qed

\begin{thm}
  Let $(V,l,r)$ be a representation of the Leibniz algebra $\g$. Then
  $$\rho=(l^*\otimes 1+ 1\otimes l)+r: \g \longrightarrow \ol(V^*\otimes V) $$is a
naive  representation  of $\g$ on  $V^*\otimes V$.
\end{thm}

\pf By Proposition \ref{pro:1cocycle}, $r:\g\longrightarrow \gl(V)$ is a $1$-cocycle on $\g$ with the coefficients in the representation $(V^*\otimes V,l^*\otimes 1+ 1\otimes l,0)$.
By Proposition \ref{pro:rep con}, $\rho=(l^*\otimes 1+ 1\otimes l)+r$ is a homomorphism from $\g$ to $\ol(V^*\otimes V)$. \qed

A {\bf trivial naive representation} $\rho_T$ of $\g$ on $\mathbb R$
is defined
 to be a homomorphism $$\rho_T=\phi+\theta:\g\longrightarrow\gl(\mathbb R)
 \oplus \mathbb R $$
 such that $\phi=0$.  By Proposition \ref{pro:rep con}, we have
\begin{pro}\label{trivial}
  Trivial naive representations of a Leibniz algebra are in one-to-one correspondence
   to $\xi\in\g^*$ such that $\xi|_{[\g,\g]_\g}=0$.
\end{pro}

The {\bf adjoint naive representation} $\adj$ is defined to be the
homomorphism $$\adj=\ad_L+\id: \g\longrightarrow\gl(\g)\oplus \g.$$

\section{Naive  cohomology of Leibniz algebras }

  Let $\rho:\g\longrightarrow\ol(V)$ be a naive representation of the Leibniz algebra $\g$. It is obvious that  $
\img(\rho) \subset\ol(V)$ is a Leibniz subalgebra so that one can
 define a
set of $k$-cochains by
$$
C^k(\g;\rho)=\{f:\otimes ^k\g\longrightarrow \img(\rho)\}$$ and
 an operator $\delta:C^k(\g;\rho)\longrightarrow
C^{k+1}(\g;\rho)$ by
\begin{eqnarray}
\nonumber\delta
c^k(x_1,\dots,x_{k+1})&=&\sum_{i=1}^k(-1)^{i+1}\Courant{\rho(x_i),c^k(x_1,\dots,\widehat{x_i},\dots,x_{k+1})}
\\\nonumber&&+(-1)^{k+1}\Courant{ c^k(x_1,\dots,x_k),\rho(x_{k+1})}\\
\label{formulapartial}&&+\sum_{1\leq i<j\leq
k+1}(-1)^ic^k(x_1,\dots,\widehat{x_i},\dots,x_{j-1},[x_i,x_j]_\frkg,x_{j+1},\dots,x_{k+1}).
\end{eqnarray}

\begin{lem}
  With the above notations, we have $\delta^2=0$.
\end{lem}

\pf For all $x\in\g$ and $u\in\img(\rho)$, define
$$
l_x(u)=\Courant{\rho(x),u},\quad r_x(u)=\Courant{u,\rho(x)}.
$$
By the fact that $\ol(V)$ is a Leibniz algebra, we can deduce that
$(\img(\rho);l,r)$ is a representation of $\g$ on $\img(\rho)$ in
the sense of Definition \ref{defi:rep old}. $\delta$ is just the
usual coboundary operator for this representation so that
$\delta^2=0.$ \qed\vspace{3mm}

Thus, we have a well-defined cochain complex
$(C^{\bullet}(\g;\rho),\delta).$  The resulting cohomology is called
the  {\bf naive  cohomology } and denoted by
$H_{naive}^\bullet(\g;\rho)$. In particular, $H_{naive}^\bullet(\g)$
and $H_{naive}^\bullet(\g;\adj)$ denote the naive cohomologies
corresponding to the trivial  naive representation and adjoint naive
representation of $\g$ respectively.

\begin{thm}
 With the above notations,   we have  $H_{naive}^\bullet(\g) =H^\bullet(\g)$.
\end{thm}
\pf If $[\g,\g]_\g=\g$, there is only one trivial naive
representation $\rho=0$
 by Proposition \ref{trivial}.
 In this case, all the cochains are also $0$. Thus, $H_{naive}^\bullet(\g)=0$. On the other hand, under the condition  $[\g,\g]_\g=\g$, it is straightforward to deduce that for any $\xi\in C^k(\g)$, $\partial \xi=0$ if and only if $\xi=0$. Thus, $H^\bullet(\g)=0$.

If $[\g,\g]_\g\neq\g$, any $0\neq\xi\in\g^*$ such that
$\xi|_{[\g,\g]_\g}=0$ gives rise to a trivial naive representation
$\rho_T$. Furthermore, we have $\img(\rho_T)=\mathbb R$ and
$C^k(\g)=\wedge^k\g^*$. Thus, the sets of cochains are the same
associated to two kinds of representations. Since $V$ is an abelian
subalgebra in $\ol(V)$, for any $\xi\in C^k(\g)$, we have
\begin{eqnarray*}
\delta
\xi(x_1,\dots,x_{k+1})&=&\sum_{i=1}^k(-1)^{i+1}\Courant{\rho_T(x_i),c^k(x_1,\dots,\widehat{x_i},\dots,x_{k+1})}
\\\nonumber&&+(-1)^{k+1}\Courant{ c^k(x_1,\dots,x_k),\rho_T(x_{k+1})}\\
\label{formulapartial}&&+\sum_{1\leq i<j\leq
k+1}(-1)^ic^k(x_1,\dots,\widehat{x_i},\dots,x_{j-1},[x_i,x_j]_\frkg,x_{j+1},\dots,x_{k+1})\\
&=&\sum_{1\leq i<j\leq
k+1}(-1)^ic^k(x_1,\dots,\widehat{x_i},\dots,x_{j-1},[x_i,x_j]_\frkg,x_{j+1},\dots,x_{k+1})\\
&=&\partial \xi(x_1,\dots,x_{k+1}).
\end{eqnarray*}
Thus, we have $H_{naive}^\bullet(\g)=H^\bullet(\g)$.\qed\vspace{3mm}

For the adjoint   naive representation $\adj=\ad_L+\id:
\g\longrightarrow\gl(\g)\oplus \g,$  any $k$-cochain $f$ is uniquely
determined by a linear map $\frkf:\otimes^k\g\longrightarrow \g$
such that
\begin{equation}\label{eq:correspondence}
f=(\ad_L\circ{\frkf},\frkf):\otimes^k\g\longrightarrow \img(\adj).
\end{equation}

\begin{thm}\label{thm:adjoint}
 With the above notations, we have $H_{naive}^\bullet(\g;\adj)=H^\bullet(\g;\ad_L,\ad_R)$.
\end{thm}
\pf Since any $k$-cochain $f:\otimes^k\g\longrightarrow \img(\adj)$ is uniquely determined
by a linear map $\frkf:\otimes^k\g\longrightarrow \g$ via \eqref{eq:correspondence}. Thus, there is a
 one-to-one correspondence between the sets of cochains associated to
 the two kinds representations via $f\leftrightsquigarrow\frkf$.
 Furthermore, we have

\begin{eqnarray*}
&&\delta
f(x_1,\dots,x_{k+1})\\&=&\sum_{i=1}^k(-1)^{i+1}\Courant{\adj(x_i),f(x_1,\dots,\widehat{x_i},\dots,x_{k+1})}
\\\nonumber&&+(-1)^{k+1}\Courant{ f(x_1,\dots,x_k),\adj(x_{k+1})}\\
\label{formulapartial}&&+\sum_{1\leq i<j\leq
k+1}(-1)^if(x_1,\dots,\widehat{x_i},\dots,x_{j-1},[x_i,x_j]_\frkg,x_{j+1},\dots,x_{k+1})\\
&=&\sum_{i=1}^k(-1)^{i+1}\Courant{\ad_L(x_i)+x_i,\ad_L\frkf(x_1,\dots,\widehat{x_i},\dots,x_{k+1})+\frkf(x_1,\dots,\widehat{x_i},\dots,x_{k+1})}
\\\nonumber&&+(-1)^{k+1}\Courant{ \ad_L\frkf(x_1,\dots,x_k)+\frkf(x_1,\dots,x_k),\ad_L(x_{k+1})+x_{k+1}}\\
\label{formulapartial}&&+\sum_{1\leq i<j\leq
k+1}(-1)^i\Big(\ad_L\frkf(x_1,\dots,\widehat{x_i},\dots,x_{j-1},[x_i,x_j]_\frkg,x_{j+1},\dots,x_{k+1})\\
&&+\frkf(x_1,\dots,\widehat{x_i},\dots,x_{j-1},[x_i,x_j]_\frkg,x_{j+1},\dots,x_{k+1})\Big)\\
&=&\sum_{i=1}^k(-1)^{i+1}\Big([\ad_L(x_i),\ad_L\frkf(x_1,\dots,\widehat{x_i},\dots,x_{k+1})]+\ad_L(x_i)\frkf(x_1,\dots,\widehat{x_i},\dots,x_{k+1})\Big)
\\\nonumber&&+(-1)^{k+1}\Big([ \ad_L\frkf(x_1,\dots,x_k),\ad_L(x_{k+1})]+\ad_L\frkf(x_1,\dots,x_k)x_{k+1}\Big)\\
\label{formulapartial}&&+\sum_{1\leq i<j\leq
k+1}(-1)^i\Big(\ad_L\frkf(x_1,\dots,\widehat{x_i},\dots,x_{j-1},[x_i,x_j]_\frkg,x_{j+1},\dots,x_{k+1})\\
&&+\frkf(x_1,\dots,\widehat{x_i},\dots,x_{j-1},[x_i,x_j]_\frkg,x_{j+1},\dots,x_{k+1})\Big)\\
&=&\sum_{i=1}^k(-1)^{i+1}\Big(\ad_L[x_i,\frkf(x_1,\dots,\widehat{x_i},\dots,x_{k+1})]_\g+[x_i,\frkf(x_1,\dots,\widehat{x_i},\dots,x_{k+1})]_\g\Big)
\\\nonumber&&+(-1)^{k+1}\Big( \ad_L[\frkf(x_1,\dots,x_k),x_{k+1}]_\g+[\frkf(x_1,\dots,x_k),x_{k+1}]_\g\Big)\\
\label{formulapartial}&&+\sum_{1\leq i<j\leq
k+1}(-1)^i\Big(\ad_L\frkf(x_1,\dots,\widehat{x_i},\dots,x_{j-1},[x_i,x_j]_\frkg,x_{j+1},\dots,x_{k+1})\\
&&+\frkf(x_1,\dots,\widehat{x_i},\dots,x_{j-1},[x_i,x_j]_\frkg,x_{j+1},\dots,x_{k+1})\Big)\\
&=&\ad_L\partial \frkf(x_1,\dots,x_{k+1})+\partial \frkf(x_1,\dots,x_{k+1}).
\end{eqnarray*}
Thus, we have $\delta f=0$ if and only if $\partial \frkf=0$. Similarly, we can prove that $f$ is exact if and only if $\frkf$ is exact. Thus, the corresponding cohomologies are isomorphic. \qed\vspace{3mm}

At last, we consider a naive representation $\rho$ such that the
image of $\rho$ is
 contained in the graph $\huaG_\varphi$ for some linear map $\varphi:V\longrightarrow\gl(V)$
 satisfying Eq. \eqref{eq:gr}.
In this case, $\rho$
 is of the form $\rho=\varphi\circ \theta+\theta$, where $\theta:\g\longrightarrow V$ is a linear map. A $k$-cochain $f:\otimes^k\g\longrightarrow \img(\rho)$ is of the form $f=\varphi\circ \frkf+\frkf$, where $\frkf:\otimes^k\g\longrightarrow V$ is a linear map.

 Define left and right actions in the sense of Definition \ref{defi:rep old} by
 \begin{eqnarray}
  \label{eq:l} l_xu&=&\pr\Courant{\rho(x),\varphi(u)+u}=\varphi(\theta(x))u;\\
  \label{eq:r}  r_xu&=&\pr\Courant{\varphi(u)+u,\rho(x)}=\varphi(u)\theta(x),
 \end{eqnarray}
 where $\pr$ is the projection from $\gl(V)\oplus V$ to $V$.  Similar to Theorem \ref{thm:adjoint}, it is easy to
 prove that
\begin{thm}
 Let $\rho$ be a naive representation such that the image of $\rho$ is contained in the graph $\huaG_\phi$ for some linear map $\varphi:V\longrightarrow\gl(V)$ satisfying Eq.  \eqref{eq:gr}. Then we have
 $$
 H_{naive}^\bullet(\g;\rho)=H^\bullet(\g;l,r),
 $$
 where $l$ and $r$ are given by \eqref{eq:l} and \eqref{eq:r} respectively.
\end{thm}

\end{document}